\documentclass[CJK,GBK,11pt]{article}
\usepackage{amsmath,amssymb,amsthm,amscd}
\usepackage[colorlinks, linkcolor=blue, anchorcolor=gree, citecolor=red]{hyperref}
\usepackage[numbers,sort&compress]{natbib}
\usepackage{graphicx}
\begin{document}
\newcommand{\Cyc}{{\rm{Cyc}}}\newcommand{\diam}{{\rm{diam}}}
\newtheorem{thm}{Theorem}[section]
\newtheorem{pro}[thm]{Proposition}
\newtheorem{lem}[thm]{Lemma}
\newtheorem{fac}[thm]{Fact}
\newtheorem{cor}[thm]{Corollary}
\newtheorem{Q}{Question}[section]
\newtheorem{C}{Conjecture}[section]
\theoremstyle{definition}
\newtheorem{definition}{Definition}[section]
\newtheorem{ex}[thm]{Example}
\newtheorem{ob}[thm]{Observtion}
\theoremstyle{remark}
\newtheorem{remark}[thm]{Remark}
\newcommand{\bth}{\begin{thm}}
\renewcommand{\eth}{\end{thm}}
\newcommand{\bex}{\begin{examp}}
\newcommand{\eex}{\end{examp}}
\newcommand{\bre}{\begin{remark}}
\newcommand{\ere}{\end{remark}}
\newcommand{\bal}{\begin{aligned}}
\newcommand{\eal}{\end{aligned}}
\newcommand{\beq}{\begin{equation}}
\newcommand{\eeq}{\end{equation}}
\newcommand{\ben}{\begin{equation*}}
\newcommand{\een}{\end{equation*}}
\newcommand{\bpf}{\begin{proof}}
\newcommand{\epf}{\end{proof}}
\renewcommand{\thefootnote}{}
\newcommand{\sdim}{{\rm sdim}}
\newcommand{\capin}[3]{\underset{#1=1}{\overset{#2}{\cap}}#3}
\def\beql#1{\begin{equation}\label{#1}}
\makeatletter
\newcommand{\rmnum}[1]{\romannumeral#1}
\newcommand{\Rmnum}[1]{\expandafter\@slowromancap\romannumeral#1@}
\makeatother
\title{\Large\bf Sylow numbers and the structure of finite groups}
\vspace{0.8cm}
\author{{\small Huaquan WEI$^{1}$\footnote{Corresponding authors. E-mail adress: weihuaquan@163.com. Projects supported by NSF of China (12061011), Guangxi (2023GXNSFAA026333, 2021GXNSFAA220105) and Guangxi Graduate Education (YCBZ2023021).}, Yi CHEN$^{1}$, Hui WU$^{1}$, Jiamin ZHANG$^{1}$, Jiawen HE$^{2}$}\\
{\tiny 1. College of Mathematics and Information Science, Guangxi University, Nanning, 530004, China}\\
{\tiny 2. College of Artificial Intelligence and Software, Nanning University, Nanning, 530200, China}\\
}
\date{}
\maketitle
\medskip
\begin{abstract}
{Suppose that the finite group $G=AB$ is a mutually permutable product of two subgroups $A$ and $B$. By using Sylow numbers of $A$ and $B$, we present some new bounds of the $p$-length $l_p(G)$ of a $p$-solvable group $G$ and the nilpotent length $F_l(G)$ and the derived length $dl(G/\Phi(G))$ of a solvable group $G$. Some known results of Zhang in J. Algebra 1995, 176 are extended.}\\

{\bf 2000 Mathematics Subject Classification}\ \ 20D10, 20D20.\\

{\bf Keywords}\ \ finite group; $p$-solvable group; solvable group; mutually permutable product; Sylow numbers

\end{abstract}


\footnote{}

\section{Introduction}\label{Sec1}
\indent

All groups considered are finite. Let $G$ be a group, we denote by $\pi(G)$ the set of prime divisors of $|G|$. Let $p\in\pi(G)$, by $n_p(G)$, we mean the Sylow $p$-number of $G$. Following \cite{Z}, we set $\tau_{p_{i}}(n)=a_i$ and $\tau(n)=max\{a_i$ $|$ $i=1,2,\ldots, t\}$ for a natural number $n$ with prime-factor decomposition $n=p_{1}^{a_{1}}p_{2}^{a_{2}}\cdots p_{t}^{a_{t}}$; for a group $G$ and a prime $p$ we set $\tau(G)=max\{\tau(n_q(G))$ $|$ $q\in\pi(G)\}$ and $\tau_p(G)=max\{\tau_p(n_q(G))$ $|$ $q\in\pi(G)\}$. For other notations and terminologies used are standard, as in \cite{G,H}.

The number of Sylow subgroups of a group $G$ is an important number pertaining to $G$. By use of Sylow numbers, many scholars have been extensively investigated the structure of finite groups. For example, a classical result due to P. Hall \cite{PH} claims that an integer $n$ is the Sylow $p$-number of a solvable group if and only if $n=p_{1}^{a_{1}}p_{2}^{a_{2}}\cdots p_{t}^{a_{t}}$ such that $p_{i}^{a_{i}}\equiv 1$ (mod $p$), $i=1,2,\ldots, t$. Zhang in \cite{Z} showed that a group $G$ is $p$-nilpotent if and only if $p$ is prime to every Sylow number of $G$ (where $p\ne 3$). Chigira in \cite{C} proved that $G$ is $3$-nilpotent if and only if $3$ is prime to every Sylow number of $G$ and $N_G(S)$ is $3$-nilpotent for every Sylow subgroup $S$ of $G$. These two results generalize some classical ones and prove affirmatively a conjecture of Huppert. Moreover, Zhang in \cite{Z} gave a bound of the $p$-length of a $p$-solvable group and a bound of the nilpotent length or derived length of a solvable group in terms of Sylow numbers. In addition, Guo and Shum in \cite{GS} also obtained a bound of the derived length of a solvable group by use of Sylow numbers.

In general, a product of two $p$-solvable (solvable) subgroups need not be $p$-solvable (solvable). However, if the group $G$ is a mutually permutable product of two $p$-solvable (solvable) subgroups, then $G$ is still a $p$-solvable (solvable) group \cite{BEA}. Recall that the product $G=AB$ of the subgroups $A$ and $B$ of a group $G$ is called a mutually permutable product of $A$ and $B$ if $AU=UA$ for any subgroup $U$ of $B$ and $BV=VB$ for any subgroup $V$ of $B$ \cite{BEA}. By use of the concept of mutually permutable product, some scholar have been studied the structure of finite groups. For instance, Cossey, Wei, Gu et al. in \cite{CL1}-\cite{GLWY} respectively obtained many bounds of the $p$-length of a mutually permutable product of two $p$-solvable groups. These results extend the celebrated Hall-Higman theorem \cite{HH} on the $p$-length of a $p$-solvable group.

In the paper, we continue the study of mutually permutable product of finite groups by use of Sylow numbers. We focus our attention on the $p$-length, nilpotent length and derived length of a mutually permutable product of two $p$-solvable or solvable groups. In detail, we present the following new bounds:

\begin{thm}\label{th1.2}
{Suppose that $G=AB$ is a mutually permutable product of two $p$-solvable subgroups $A$ and $B$, where $p\in\pi(G)$. If one of the following conditions holds, then the $p$-length $l_p(G)\leq max\{1+\frac{\tau_p(A)}{2}, 1+\frac{\tau_p(B)}{2}\}:$
	
$(a)$ $(|G|, p-1)=1$;

$(b)$ either $A$ or $B$ is $p$-nilpotent.
}
\end{thm}

For simplicity, we define the following real functions:
	\begin{equation*}
	f(x)=
	\begin{cases}
	\log_{3}(x/2) &x>0,\\
	\log_{3}(1/2) &x=0.
	\end{cases}	
	\end{equation*}	

	\begin{equation*}
	g(x)=
	\begin{cases}
	\log_{2}(x) &x>0,\\
	0 &x=0.
	\end{cases}	
	\end{equation*}	

\begin{thm}\label{th1.3}
{Suppose that $G=AB$ is a mutually permutable product of two solvable subgroups $A$ and $B$. Then

$(a)$ the nilpotent length $F_l(G)\leq max\{4+2f(\tau(A)), 4+2f(\tau(B))\};$

$(b)$ the derived length $dl(G/\Phi(G))\leq max\{2+6g(\tau(A)), 2+6g(\tau(B))\}$.
}
\end{thm}

\begin{remark}\label{r1.4}
{(1) The equality can be satisfied in Theorem \ref{th1.2} for some proper subgroups $A$ and $B$. For example, let $G=S_4$. Then $G=AB$ is a mutually permutable product of $A$ and $B$, where $A\cong A_4$, $B\cong D_8$. Clearly, $\tau_2(A)=2$, $\tau_2(B)=0$ and $l_2(G)=2$, hence $l_2(G)=1+\frac{\tau_p(A)}{2}$.

(2) There exists a solvable group $G$ such that $F_l(G)=dl(G)=2$. For example, $G=S_3=C_3\rtimes C_2$ with $p=2$. Hence the bounds in Theorem \ref{th1.3} cannot be improved respectively to be $max\{3+2f(\tau(A)), 3+2f(\tau(B))\}$ and $max\{1+6g(\tau(A)), 1+6g(\tau(B))\}$.}

(3) We don't know whether the condition (a) or (b) in Theorem \ref{th1.2} can be removed. So we pose the following conjecture.

\end{remark}

\begin{C}\label{C1.4}
{Suppose that $G=AB$ is a mutually permutable product of two $p$-solvable subgroups $A$ and $B$, where $p\in\pi(G)$. Then $l_p(G)\leq max\{1+\frac{\tau_p(A)}{2}, 1+\frac{\tau_p(B)}{2}\}$.
}
\end{C}

\medskip
\section{Preliminaries}\label{Sec2}
\indent

Let $G$ be a group and let $\pi$ be a set of primes. It is well-known that $O^{\pi}(G)$ is the intersection of all normal subgroups $N$ of $G$ such that $G/N$ is a $\pi$-group. Hence $G/O^{\pi}(G)$ is the maximal $\pi$-factor group of $G$ (\cite[IX, 1.1]{HB}). Following \cite{CL2}, we invoke the following definition way of $p$-length of a $p$-solvable group.

If $p$ is a prime, the lower $p$-series of $G$ is $$G\ge O^{p'}(G)\ge O^{p',p}(G)\ge O^{p',p,p'}(G)\ge\cdots.$$If $G$ is $p$-solvable, the last term of the lower $p$-series is 1 and if the lower $p$-series of $G$ is $$G=G_0\ge G_1\ge\cdots\ge G_s=1,$$then the $p$-length of $G$ is the number of non-trivial $p$-groups in the set $$\{G/G_1, G_1/G_2, \ldots, G_{s-1}/G_s\}.$$

Now let $G$ be a solvable group. So-called the derived length $dl(G)$ and the nilpotent length $F_l(G)$ of $G$ are the length of a shortest abelian series and the length of the upper nilpotent series in $G$, respectively.

\begin{lem}\label{L1} $($\cite[Theorem 4.1.15]{BEA}$)$
{Let the group $G$ be the product of the mutually permutable subgroups $A$ and $B$. If $A$ and $B$ are $p$-solvable, then $G$ is $p$-solvable.}
\end{lem}

\begin{lem}\label{L2} $($\cite[Lemma 4.1.10]{BEA}$)$
{Let the group $G$ be the product of the mutually permutable subgroups $A$ and $B$. If $N$ is a normal subgroup of $G$, then $G/N$ is a mutually permutable product of $AN/N$ and $BN/N$.}
\end{lem}

\begin{lem}\label{L3} $($\cite[Theorem 4.3.11]{BEA}$)$
{Let the non-trivial group $G$ be the product of mutually permutable subgroups $A$ and $B$. Then $A_GB_G$ is not trivial.}
\end{lem}

\begin{lem}\label{L4} $($\cite[Lemma 4.3.3]{BEA}$)$
{Let the group $G$ be the product of the mutually permutable subgroups $A$ and $B$. Then
	
	$(1)$ If $N$ is a minimal normal subgroup of $G$, then $\{A\cap N, B\cap N\}\subseteq\{N, 1\}$.
	
	$(2)$ If $N$ is a minimal normal subgroup of $G$ contained in $A$ and $B\cap N=1$, then $N\le C_G(A)$ or $N\le C_G(B)$. If furthermore $N$ is not cyclic, then $N\le C_G(B)$.}
\end{lem}

\begin{lem}\label{L5} $($\cite[Corollary 4.1.22]{BEA}$)$
{Let the group $G$ be the product of the mutually permutable subgroups $A$ and $B$. Then

	$(1)$ If $U$ is a subgroup of $G$, then $(A\cap U)(B\cap U)$ is the mutually permutable product of $A\cap U$ and $B\cap U$;

	$(2)$ If $U$ is a normal subgroup of $G$, then $(A\cap U)(B\cap U)$ is a normal subgroup of $G$.}
\end{lem}

\begin{lem}\label{L6}
{Let $A$ and $H$ be subgroups of a group $G$ such that $G=AH$, $A\cap H=1$. If $A$ is a normal $p$-group of $G$, then $|H:N_H(Q)|=|G:AN_G(Q)|$ for any $Q\in$ Syl$_q(H) (q\neq p)$.}
\end{lem}

\proof Let $G$ be a counter-example of minimal order. Note that $G_1=AN_G(Q)=A(H\cap G_{1})$ satisfies the hypotheses. So we consider the following two cases:

    \textbf{Case 1.} $|G_1|<|G|$. By the minimality of $G$, we have $$|H\cap G_1:N_{H\cap G_1}(Q)|=|G_1:AN_{G_{1}}(Q)|=1.$$ We observe that $N_{H\cap G_1}(Q)=H\cap N_G(Q)$, this indicates that $$H\cap AN_G(Q)=H\cap N_G(Q)=N_H(Q).$$ Noticing that $AN_G(Q)=A(H\cap AN_G(Q))$, hence $$|AN_G(Q)/A|=|H\cap AN_G(Q)|=|N_H(Q)|$$ and so $$|G:AN_G(Q)|=|G/A:AN_G(Q)/A|=|H:N_H(Q)|,$$which is contrary to the choice of $G$.

    \textbf{Case 2.} $G=G_1=AN_G(Q)$. In this case, $$N_G(Q)/(A\cap N_G(Q))\cong G/A\cong H.$$ By the uniqueness of the Sylow $q$-subgroup of $N_G(Q)/(A\cap N_G(Q))$, we obtain $Q\unlhd H$ and hence $|H:N_H(Q)|=1=|G:AN_G(Q)|$, a contradiction.\qed

\begin{lem}\label{L7} $($\cite[Theorem 10]{Z}$)$ Let $G$ be a $p$-solvable group. Then the $p$-length $l_p(G)\leq 1+\frac{\tau_p(G)}{2}$.

\end{lem}

\medskip
\section{Proofs of main results}\label{Sec3}

Proof of Theorem \ref{th1.2}:

Let $G$ be a counter-example of minimal order. We proceed in steps.

    $(1)$ $G$ is $p$-solvable.

    This follows from Lemma \ref{L1}.

    $(2)$ $N=O_p(G)=C_G(N)$ is unique minimal normal and complemented in $G$.

    Let $N$ be a minimal normal subgroup of $G$. We consider $\overline{G}=G/N$ together with $\overline{A}=AN/N$ and $\overline{B}=BN/N$. By Lemma \ref{L2}, $\overline{G}$ is the mutually product of two $p$-solvable subgroups $\overline{A}$ and $\overline{B}$. It is clear that $(|\overline{G}|, p-1)=1$ if $G$ satisfies (a) and either $\overline{A}$ or $\overline{B}$ is $p$-nilpotent if $G$ satisfies (b). Hence $\overline{G}$ satisfies the hypotheses of Theorem \ref{th1.2}. The choice of $G$ and \cite[Lemma 1]{Z} imply that $$l_p(\overline{G})\le max\{1+\tau_p(\overline{A})/2, 1+\tau_p(\overline{B})/2\}\le max\{1+\tau_p(A)/2, 1+\tau_p(B)/2\}.$$If $N_1$ is minimal normal in $G$ with $N_1\ne N$, then we also have $$l_p(G/N_1)\le max\{1+\tau_p(A)/2, 1+\tau_p(B)/2\}.$$It follows that $$l_p(G)\le max\{l_p(G/N), l_p(G/N_1)\}\le max\{1+\tau_p(A)/2, 1+\tau_p(B)/2\},$$
    which is contrary to the choice of $G$. Hence $N$ is the unique minimal normal subgroup of $G$. Moreover, if $N\le O_{p'}(G)$ or $N\le\Phi(G)$, then $$l_p(G)=l_p(\overline{G})\le max\{1+\tau_p(A)/2, 1+\tau_p(B)/2\},$$a contradiction. Hence $O_{p'}(G)=\Phi(G)=1$ and $N=O_p(G)=C_G(N)$ is complemented in $G$, as desired.

    $(3)$ $N\le A\cap B$.

    Since $A_GB_G\ne 1$ by Lemma \ref{L3}, we may assume $N\le A$ by $(2)$. If $N\not\le B$, then $N\cap B=1$ by Lemma \ref{L4}(1). If $N$ is cyclic, then $N=C_G(N)\in$ Syl$_p(G)$, hence $l_p(G)=1$, a contradiction. Thus $N$ is not cyclic and $N\le C_G(B)$ by Lemma \ref{L4}(2). Furthermore, $B\le C_G(N)=N\le A$ and thereby $G=AB=A$. In this case, Theorem \ref{th1.2} is true by Lemma \ref{L7}, also a contradiction. This proves $N\le A\cap B$.

    $(4)$ $O^{p'}(G)=G$.

    If not, then $O^{p'}(G)<G$. In view of Lemma \ref{L5}, $(A\cap O^{p'}(G))(B\cap O^{p'}(G))$ is normal in $G$ and it is the mutually permutable product of $A\cap O^{p'}(G)$ and $B\cap O^{p'}(G)$. Clearly, $O^{p'}(A)\le A\cap O^{p'}(G)$ and $O^{p'}(B)\le B\cap O^{p'}(G)$. Hence $O^{p'}(G)=(A\cap O^{p'}(G))(B\cap O^{p'}(G))$ satisfies the hypotheses of Theorem \ref{th1.2}. The minimality of $G$  and \cite[Lemma 1]{Z} imply that

    \begin{eqnarray*}
    l_p(G)&=&l_p(O^{p'}(G))\le max\{1+\dfrac{\tau_p(A\cap O^{p'}(G))}{2}, 1+\dfrac{\tau_p(B\cap O^{p'}(G))}{2}\}
    \\
    &\le& max\{1+\dfrac{\tau_p(A)}{2}, 1+\dfrac{\tau_p(B)}{2}\},
    \end{eqnarray*}
    a contradiction.

    $(5)$ $G$ has a maximal subgroup $L$ such that $N\cap L=1$ and $G=NL$. Moreover, $\tau_p(L_1)\leq \tau_p(A)-2$ and $\tau_p(L_2)\leq \tau_p(B)-2$, where $L_1=A\cap L$ and $L_2=B\cap L$.

    Such a subgroup $L$ of $G$ does exist by (2). If both $A$ and $B$ are $p$-nilpotent, then both $A$ and $B$ are $p$-groups by $N=C_{G}(N)$; as a result, $G=AB$ is a $p$-group and hence $l_p(G)=1$, which is impossible. Now we prove $\tau_p(L_1)\leq \tau_p(A)-2$. Suppose otherwise $\tau_p(L_1)\geq \tau_p(A)-1$. We consider the following three cases:

    \textbf{Case 1.} $A$ is not $p$-nilpotent and $B$ is $p$-nilpotent. In this case, $B$ is a $p$-group. Since $A=NL_1$ and $N\cap L_1=1$, there exists $Q\in$ Syl$_q(L_1)(q\neq p)$ such that $$\tau_p(|L_1:N_{L_{1}}(Q)|)\geq \tau_p(|A:N_{A}(Q)|)-1.$$ By Lemma \ref{L6}, we have $$|L_1:N_{L_{1}}(Q)|=|A:NN_{A}(Q)|=\dfrac{|A:N_{A}(Q)|}{|N:N_N(Q)|}.$$ It follows that $\tau_p(|L_1:N_{L_{1}}(Q)|)= \tau_p(|A:N_{A}(Q)|)-1$ and $|N/N_N(Q)|=p$. Write $D=N_N(Q)$. Then $D=C_N(Q)$ and $Q$ is obviously bound to act faithfully on $N/D$. Thus $Q$ is cyclic and $q$ $|$ $p-1$. This is contrary to $(|G|,p-1)=1$ if $G$ satisfies (a). Now assume that $G$ satisfies (b). We claim that $N_{L_1}(Q)=C_{L_1}(Q)$. In fact, if $N_{L_1}(Q)>C_{L_1}(Q)$, then $Q\leq (N_{L_1}(Q))'$; it follows that $Q$ acts trivially on $N/D$, which is impossible. Hence $N_{L_1}(Q)=C_{L_1}(Q)$ and $L_1$ is $q$-nilpotent, of course, $A$ is also $q$-nilpotent. Since $q<p$ and $Q$ is cyclic, $QB$ is $q$-nilpotent, namely $B$ is normalized by $Q$. It follows that $G$ is $q$-nilpotent, which is contrary to $O^{p'}(G)=G$.

    \textbf{Case 2.}  $A$ is not $p$-nilpotent and $B$ is not $p$-nilpotent. In the present case, $G$ must not satisfy (b), so $G$ satisfies (a) by the hypotheses of Theorem \ref{th1.2}. With similar arguments as in Case 1, we can get a contradiction.

    \textbf{Case 3.}  $A$ is $p$-nilpotent and $B$ is not $p$-nilpotent. In this case, $A$ is a $p$-group. If $\tau_p(L_2)\leq \tau_p(B)-2$ then by the minimality of $G$, we have
\begin{eqnarray*}
    l_p(G)&=& 1+l_p(G/N)\le 1+max\{1+\dfrac{\tau_p(A/N)}{2}, 1+\dfrac{\tau_p(B/N)}{2}\}
    \\&=&1+max\{1,1+\dfrac{\tau_p(L_2)}{2}\}=2+\dfrac{\tau_p(L_2)}{2}\le 1+\dfrac{\tau_p(B)}{2}
    \\&=& max\{1+\dfrac{\tau_p(A)}{2},1+\dfrac{\tau_p(B)}{2}\},
    \end{eqnarray*}
    a contradiction. Now assume $\tau_p(L_2)\geq \tau_p(B)-1$. With similar arguments as in Case 1, we can also derive a contradiction.

    $(6)$ Finishing the proof.

    By $(5)$, we obtain

    \begin{eqnarray*}
    l_p(G)&=& 1+l_p(G/N)\le 1+max\{1+\dfrac{\tau_p(A/N)}{2}, 1+\dfrac{\tau_p(B/N)}{2}\}
    \\&=& 1+max\{1+\dfrac{\tau_p(L_1)}{2}, 1+\dfrac{\tau_p(L_2)}{2}\}
    \\&\le& max\{1+\dfrac{\tau_p(A)}{2},1+\dfrac{\tau_p(B)}{2}\}.
    \end{eqnarray*}

    This is the final contradiction and the proof is complete.\qed

Proof of Theorem \ref{th1.3}:

Let $G$ be a counter-example of minimal order. We proceed in steps.

    $(1)$ $G$ is solvable.

    This follows from Lemma \ref{L1}.

    $(2)$ $N=O_p(G)=F(G)$ is unique minimal normal and complemented in $G$ for some $p\in\pi(G)$ and $N=C_G(N)$, $\Phi(G)=1$.

    $(a)$ Let $N$ be a minimal normal subgroup of $G$. We consider $\overline{G}=G/N$ together with $\overline{A}=AN/N$ and $\overline{B}=BN/N$.  By Lemma \ref{L2}, $\overline{G}$ is the mutually product of two solvable subgroups $\overline{A}$ and $\overline{B}$, hence $\overline{G}$ satisfies the hypotheses of Theorem \ref{th1.3}. The choice of $G$ and \cite[Lemma 1]{Z} imply that $$F_l(\overline{G})\le max\{4+2f(\tau(\overline A)), 4+2f(\tau(\overline B))\}\le max\{4+2f(\tau(A)), 4+2f(\tau(B))\}.$$If $N_1$ is minimal normal in $G$ with $N_1\ne N$, then we also have $$F_l(G/N_1)\le max\{4+2f(\tau(A)), 4+2f(\tau(B))\}.$$It follows that $$F_l(G)\le max\{F_l(G/N), F_l(G/N_1)\}\le max\{4+2f(\tau(A)), 4+2f(\tau(B))\},$$
    a contradiction. Therefore $N$ is the unique minimal normal subgroup of $G$. In view of $(1)$, $N\le O_{p}(G)$ for some $p\in\pi(G)$. Moreover, if $N\le\Phi(G)$, then $F(G/N)=F(G)/N>1$ and $$F_l(G)=F_l(\overline{G})\le max\{4+2f(\tau(A)), 4+2f(\tau(B))\},$$contradicting to the choice of $G$. Hence $\Phi(G)=1$ and $N=O_p(G)=F(G)=C_G(N)$ is complemented in $G$.

    $(b)$ Let $N$ be a minimal normal subgroup of $G$. If $\Phi(G)>1$, then $G/\Phi(G)$ satisfies the hypotheses of Theorem \ref{th1.3}. By \cite[Lemma 1]{Z}, we have
    \begin{eqnarray*}
    dl(G/\Phi(G))&=&dl(\frac{G/\Phi(G)}{\Phi(G/\Phi(G))})\\&\le& max\{2+6g(\tau(\overline{A})),2+6g(\tau(\overline{B}))\}\\
    &\le& max\{2+6g(\tau(A)),2+6g(\tau(B))\},
    \end{eqnarray*}
    which is contrary to the choice of $G$. Hence $\Phi(G)=1$. Let $U$ be the set of maximal subgroups of $G$ and let $N^{\Rmnum{1}}$ be the set of maximal subgroups of $G$ containing $N$. Write $N^{\Rmnum{2}}=U-N^{\Rmnum{1}}$. It is clear that both $N^{\Rmnum{1}}$ and $N^{\Rmnum{2}}$ are not empty and $L=\underset{M\in N^{\Rmnum{2}}}{\bigcap}M\unlhd G$ by $N\unlhd G$. If $L>1$, then $G$ has a minimal normal subgroup $R$ such that $R\leq L$. Similarly, we denote $R^{\Rmnum{2}}=U-R^{\Rmnum{1}}$, where $R^{\Rmnum{1}}$ is the set of maximal subgroups of $G$ containing $R$. Noth that $$(\underset{M\in N^{\Rmnum{1}}}{\bigcap}M)\cap(\underset{M\in R^{\Rmnum{1}}}{\bigcap}M)\subseteq(\underset{M\in N^{\Rmnum{1}}}{\bigcap}M)\cap(\underset{M\in N^{\Rmnum{2}}}{\bigcap}M)=\Phi(G)=1.$$ So if we write $\Phi(G/N)=L_1/N$ and $\Phi(G/R)=L_2/R$, then $L_1\cap L_2=1$ and
    \begin{eqnarray*}
    dl(G) &\leq& max\{dl(G/L_1),dl(G/L_2)\}=max\{dl(\dfrac{G/N}{\Phi(G/N)}),dl(\dfrac{G/R}{\Phi(G/R)})\}\\
    &\leq& max\{2+6g(\tau(A)),2+6g(\tau(B))\},
    \end{eqnarray*}
    a contradiction. Hence $L=1$. Take $M_0\in N^{\Rmnum{2}}$. Then $G=NM_0$ and $N\cap M_0 =1.$ If $(M_0)_G=1$, then $C_G(N)\cap M_0\leq (M_0)_G=1$. Thereby $N=C_G(N)$ by $N\le C_G(N)$. If $(M_0)_G>1$, then we can take a minimal normal subgroup $N_1$ of $G$ contained in $(M_0)_G$. Since $$N^{\Rmnum{2}}=(N_{1}^{\Rmnum{1}}+N_{1}^{\Rmnum{2}})\cap N^{\Rmnum{2}}=(N_{1}^{\Rmnum{1}}\cap N^{\Rmnum{2}})+(N_{1}^{\Rmnum{2}}\cap N^{\Rmnum{2}}),$$ either $N^{\Rmnum{2}}\subseteq N_{1}^{\Rmnum{2}}$ or $N^{\Rmnum{2}}\subseteq N_{1}^{\Rmnum{1}}$ if either $N_{1}^{\Rmnum{1}}\cap N^{\Rmnum{2}}$ is empty or $N_{1}^{\Rmnum{2}}\cap N^{\Rmnum{2}}$ is empty. This would result either $1<N_1\le M_0\in N_{1}^{\Rmnum{2}}$ or $N_1\leq \underset{M\in N_{1}^{\Rmnum{1}}}{\bigcap}M\leq \underset{M\in N^{\Rmnum{2}}}{\bigcap}M=1$, which is absurd. So both $N_{1}^{\Rmnum{1}}\cap N^{\Rmnum{2}}$ and $N_{1}^{\Rmnum{2}}\cap N^{\Rmnum{2}}$ are not empty. Obviously $\underset{M\in(N^{\Rmnum{2}}\cap N_{1}^{\Rmnum{2}}) }{\bigcap}M\unlhd G$ by $N\unlhd G$ and $N_1\unlhd G$. If $\underset{M\in(N^{\Rmnum{2}}\cap N_{1}^{\Rmnum{2}}) }{\bigcap}M>1$, then we can also take a minimal normal subgroup $R_1$ of $G$ such that $R_1\leq \underset{M\in(N^{\Rmnum{2}}\cap N_{1}^{\Rmnum{2}}) }{\bigcap}M$. Write $\Phi(G/N_1)=L_{11}/N_1$ and $\Phi(G/R_1)=L_{12}/R_1$. Since
    \begin{eqnarray*}
    (\underset{M\in N_{1}^{\Rmnum{1}} }{\bigcap}M)\cap(\underset{M\in R_{1}^{\Rmnum{1}} }{\bigcap}M) &\subseteq&  (\underset{M\in (N_{1}^{\Rmnum{1}}\cap N^{\Rmnum{2}})}{\bigcap}M)\cap(\underset{M\in (N_{1}^{\Rmnum{2}}\cap N^{\Rmnum{2}})}{\bigcap}M)=\underset{M\in N^{\Rmnum{2}}}{\bigcap}M=1,
    \end{eqnarray*}
    $L_{11}\cap L_{12}=1$. This implies that
    \begin{eqnarray*}
    dl(G) &\leq& max\{dl(G/L_{11}),dl(G/L_{12})\}=max\{dl(\dfrac{G/N_1}{\Phi(G/N_1)}),dl(\dfrac{G/R_1}{\Phi(G/R_1)})\}\\
    &\leq& max\{2+6g(\tau(A)),2+6g(\tau(B))\},
    \end{eqnarray*}
    also a contradiction. Hence $\underset{M\in(N^{\Rmnum{2}}\cap N_{1}^{\Rmnum{2}})}{\bigcap}M=1$. Noticing that $N_{1}^{\Rmnum{2}}\cap N^{\Rmnum{2}}\ne N^{\Rmnum{2}}$ and $N_{1}^{\Rmnum{2}}\cap N^{\Rmnum{2}}$ is not empty, so by repeating the above steps, we see that there exists a positive integer $\alpha$ such that $$(M_{\alpha})_G=\underset{M\in(N^{\Rmnum{2}}\cap N_{1}^{\Rmnum{2}}\cap\dots\cap N_{\alpha}^{\Rmnum{2}})}{\bigcap}M=1,$$ and $M_\alpha\in N^{\Rmnum{2}}\cap N_{1}^{\Rmnum{2}}\cap\dots\cap N_{\alpha}^{\Rmnum{2}}$. Hence $G=NM_\alpha$ and $N\cap M_\alpha=1$. Obviously, $C_G(N)\cap M_\alpha\leq (M_\alpha)_G=1$, so $N=C_G(N)$ by $N\le C_G(N)$. Thus $N=O_p(G)=F(G)$ is the unique minimal normal subgroup of $G$ and is complemented in $G$.

    $(3)$ Write $|N|=p^{n}$, where $n$ is a positive integer. Then $n>1$.

    By $(2)$, there exists a maximal subgroup $K$ of $G$ such that $N\cap K=1$ and $G=NK$. Noth that $K$ acts faithfully on space $N$ by $N=C_{G}(N)$. Hence $K$ is a linear group of degree $n$. But $N$ is a minimal normal subgroup of $G=NK$, so $K$ is a finite solvable completely reducible linear group of degree $n$. If $n=1$, then $K$ is an abelian group. Therefore $$F_l(G)=F_l(K)+1=2\leq max\{4+2f(\tau(A)), 4+2f(\tau(B))\}$$ and $$dl(G)\leq dl(K)+1=2\leq max\{2+6g(\tau(A)), 2+6g(\tau(B))\},$$ a contradiction. Hence $n>1$.

    $(4)$ $N\le A\cap B$.

    Since $A_GB_G\ne 1$ by Lemma \ref{L3}, we may assume $N\le A$ by $(2)$. If $N\not\le B$, then $N\cap B=1$ by Lemma \ref{L4}(1). Again, $N$ is not cyclic by $(3)$, so $N\le C_G(B)$ by Lemma \ref{L4}(2). Furthermore, $B\le C_G(N)=N\le A$ and so $G=AB=A$. By $(2)$, there exists a maximal subgroup $K$ of $G$ such that $N\cap K=1$ and $G=NK$. If $F(K)$ is a $p$-group, then $NF(K)\unlhd G$ is also a $p$-group. Hence $NF(K)=N$ and $F(K)\leq N\cap K=1$. This contradicts the fact that $K$ is a nontrivial solvable group. So $N=C_N(O_q(K))\times[N,O_q(K)]$ for some $p\ne q\in\pi(F(K))$. Because $C_N(O_q(K))$ is normalized by $N$ and $K$, we have $C_N(O_q(K))\unlhd G$. Noticing that $[N,O_q(K)]\neq 1$ by $(2)$ and $N$ is minimal normal in $G$, thereby $C_N(O_q(K))=1$. It follows that $N\cap N_G(Q)=1$ where $Q\in$ Syl$_q(K)$ and $$|G:N_G(Q)|=|N|\cdot\dfrac{|G|}{|NN_G(Q)|}=|N||G:NN_G(Q)|.$$ This proves that $|N|$ divides $|G:N_G(Q)|$ and $n\leq\tau(G)=\tau(A)$. Noth that $K$ is a finite solvable completely reducible linear group of degree $n$. In the following, we respectively consider the nilpotent length $F_l(G)$ and the derived length $dl(G/\Phi(G))$:

    $(a)$ If $F_{l}(K)=1$, then $$F_l(G)=F_l(K)+1=2\leq max\{4+2f(\tau(A)), 4+2f(\tau(B))\},$$ a contradiction. So $F_l(K)>1$. By $(2)$, $n\le\tau(A)$ and \cite{HT}, we have
    \begin{eqnarray*}
    F_l(G)-1&=&F_l(K)\le 3+2\log_{3}(n/2)=3+2f(n)\le 3+2f(\tau(A))\\
    &\le& max\{3+2f(\tau(A)), 3+2f(\tau(B))\},
    \end{eqnarray*}
    again a contradiction.

    $(b)$ By $(2)$, $(3)$, $n\le\tau(A)$ and \cite{HB1957}, we obtain
    \begin{eqnarray*}
    dl(G/\Phi(G))-1&=&dl(G)-1\le dl(K)\le 6\log_{2}(n)=6g(n)\le 6g(\tau(A))\\
    &\le& max\{6g(\tau(A)), 6g(\tau(B))\},
    \end{eqnarray*}
    a contradiction. This shows $N\le A\cap B$.

    $(5)$ Either $n\leq\tau(A)$ or $n\leq\tau(B)$.

    If $G=N$, then $F_l(G)=dl(G/\Phi(G))=1$ by $N$ is an elementary abelian $p$-group. Noth that $$max\{4+2f(\tau(A)), 4+2f(\tau(B))\}\geq 4+2f(0)=4+2\log_{3}(1/2)>2$$ and $$max\{2+6g(\tau(A)), 2+6g(\tau(B))\}\geq 2+6g(0)=2.$$ This is contrary to the choice of $G$. So $G>N$ and $(A/N)_{G/N}(B/N)_{G/N}\neq 1$ by Lemma \ref{L3}. Without loss of generality, we may assume $(A/N)_{G/N}\neq 1$. Then $N<A_G$, $A_G=NK_1$ and $N\cap K_1=1$, where $K_1=K\cap A_G\ne 1$. If $O_p(K_1)=F(K_1)>1$, then $NF(K_1)$ is a normal $p$-subgroup of $G$,  so $F(K_1)\leq N\cap K_1=1$ by $(2)$, a contradiction. Therefore, there exists $1\neq O_q(K_1)\leq F(K_1)$ where $p\neq q\in\pi(K_1)$ such that $N=C_N(O_q(K_1))\times [N,O_q(K_1)]$. Noticing that $O_q(K_1)$ $char$ $K_1\unlhd K$ and $N$ is abelian, hence $C_N(O_q(K_1))\unlhd G$ and $C_N(O_q(K_1))=1$ or $N$ by $(2)$. If $C_N(O_q(K_1))=N$, then $O_q(K_1)\le C_G(N)=N$, which is impossible. Hence $C_N(O_q(K_1))=1$. Let $Q\in$ Syl$_q(K_1)$. Then $Q\in$ Syl$_q(A_G)$, $O_q(K_1)\leq Q$ and $N\cap N_{A_G}(Q)=C_N(Q)\leq C_N(O_q(K_1))=1$, so $$|A_G:N_{A_G}(Q)|=|N|\cdot\dfrac{|A_G|}{|NN_{A_G}(Q)|}=|N||A_G:NN_{A_G}(Q)|.$$ Hence $n\leq\tau(A_G)\leq\tau(A)$, as desired.

    $(6)$ Finishing the proof.

    Without loss of generality, we may assume $n\leq\tau(A)$ by $(5)$. By $(2)$, there exists a maximal subgroup $K$ of $G$ such that $N\cap K=1$ and $G=NK$. Noth that $K$ is a finite soluble completely reducible linear group of degree $n$.

    $(a)$ If $F_{l}(K)=1$, then $$F_l(G)=F_l(K)+1=2\leq max\{4+2f(\tau(A)), 4+2f(\tau(B))\},$$ a contradiction. So $F_l(K)>1$ and by $(2)$, $n\le\tau(A)$ and \cite{HT}, we have
    \begin{eqnarray*}
    F_l(G)-1&=&F_l(K)\le 3+2\log_{3}(n/2)=3+2f(n)\le 3+2f(\tau(A))\\
    &\le& max\{3+2f(\tau(A)), 3+2f(\tau(B))\},
    \end{eqnarray*}
    also a contradiction.

    $(b)$ By $(2)$, $(3)$, $n\le\tau(A)$ and \cite{HB1957}, we obtain
    \begin{eqnarray*}
    dl(G/\Phi(G))-1&=&dl(G)-1\le dl(K)\le 6\log_{2}(n)=6g(n)\le 6g(\tau(A))\\
    &\le& max\{6g(\tau(A)), 6g(\tau(B))\}.
    \end{eqnarray*}
    This is the final contradiction and the proof is complete.\qed
\noindent


\end{document}